\documentclass[oneside,12pt]{amsart}
\usepackage{amscd,amsmath,amssymb,euscript}
\usepackage[frame,cmtip,curve,arrow,matrix,line,graph]{xy}

\title{Spaces of stability conditions}
\author{Tom Bridgeland}

\address{Department of Pure Mathematics, University of Sheffield, Hicks Building,
Hounsfield Road, Sheffield, S3 7RH, UK.}

\email{t.bridgeland@sheffield.ac.uk}

\newtheorem{thm}{Theorem}[section]

\newtheorem{prop}[thm]{Proposition}
\newtheorem{lemma}[thm]{Lemma}
\newenvironment{pf}{\paragraph{Proof}}{\qed\par\medskip}
\theoremstyle{definition}
\newtheorem{defn}[thm]{Definition}
\newtheorem{remark}[thm]{Remark}

\newtheorem{example}[thm]{Example}

\newcommand{\Pol}{\operatorname{\mathcal{P}}}
\newcommand{\NS}{\operatorname{NS}}

\newcommand{\GL}{\operatorname{GL}}
\newcommand{\GLp}{\operatorname{GL^+}}
\newcommand{\grp}{{\tilde{\GLp}}(2,\R)}
\newcommand{\Con}{\operatorname{Conf}}

\renewcommand{\leq}{\leqslant}
\renewcommand{\geq}{\geqslant}
\newcommand{\ch}{\operatorname{ch}}
\newcommand{\End}{\operatorname{End}}

\newcommand{\D}{\operatorname{\mathcal{D}}}
\newcommand{\Coh}{\operatorname{Coh}}
\newcommand{\Aut}{\operatorname{Aut}}
\newcommand{\isom}{\cong}
\newcommand{\PSL}{\operatorname{PSL}}
\newcommand{\h}{\mathfrak{h}}
\newcommand{\g}{\mathfrak{g}}
\newcommand{\hreg}{\h^{\text{reg}}}
\newcommand{\Mod}{\operatorname{Mod}}
\newcommand{\Def}{\operatorname{Def}}

\newcommand{\tensor}{\otimes}
\newcommand{\PP}{\operatorname{\mathbb P}}
\newcommand{\C}{\mathbb C}
\newcommand{\T}{\mathcal T}

\newcommand{\F}{\mathcal F}

\newcommand{\Z}{\mathbb Z}
\newcommand{\A}{\mathcal A}

\newcommand{\OO}{\mathcal O}

\newcommand{\Hom}{\operatorname{Hom}}
\newcommand{\eu}{\operatorname{\chi}}

\newcommand{\LL}{\mathcal{L}}

\newcommand{\lra}{\longrightarrow}
\newcommand{\R}{{\mathbb{R}}}
\renewcommand{\L}{{\operatorname{L}}}

\newcommand{\Stab}{\operatorname{Stab}}
\newcommand{\td}{\operatorname{td}}

\renewcommand{\Im}{\operatorname{Im}}
\renewcommand{\Re}{\operatorname{Re}}
\renewcommand{\P}{\mathcal{P}}
\newcommand{\ZZ}{\mathcal{Z}}

\newcommand{\SL}{\operatorname{SL}}
\newcommand{\Q}{\mathcal{Q}}

\newcommand{\M}{\operatorname{\mathcal{M}}}
\newcommand{\N}{\operatorname{\mathcal{N}}}
\newcommand{\QQ}{\mathbb{Q}}
\newcommand{\RR}{\operatorname{R}}

\newcommand{\Fuk}{\operatorname{Fuk}}
\newcommand{\vol}{\operatorname{vol}}
\newcommand{\HH}{\mathcal{H}}

\newcommand{\Gammah}{\hat{\Gamma}}
\newcommand{\Dh}{\hat{\D}}

\newcommand{\hh}{\hat{\h}}
\newcommand{\hregh}{\hat{\h}^{\text{reg}}}

\newcommand{\Wh}{\hat{W}}

\begin{document}

\begin{abstract}
Stability conditions are a mathematical way to understand $\Pi$-stability for D-branes in string theory.
Spaces of stability conditions seem to be related to moduli spaces of conformal field theories.
This is a survey article describing what is currently known about spaces of stability conditions, and giving some pointers for future research.
\end{abstract}

\maketitle

\section{Introduction}

Stability conditions on triangulated categories were introduced in  \cite{Br1}. The motivation was to understand Douglas' work on $\Pi$-stability
for D-branes in string theory \cite{Do2,Do3}. Since then a fair number of examples have been computed  and the definition has been further scrutinized.
The aim of this paper is twofold: firstly to survey the known examples of spaces of stability conditions, and secondly to
float the idea that there is some yet-to-be-discovered  construction that will allow one to define interesting geometric structures on these spaces. Understanding this construction seems to me to be the logical next step, but there are certainly other interesting questions that can be asked, and throughout the paper I have tried to point out some of the many parts of the story that have yet to be properly understood.

It is perhaps fair to say that the whole subject of stability conditions has a slightly temporary feel to it because the definition itself looks a bit unnatural.
On the other hand many of the examples that have been computed are extremely neat, and the whole idea of extracting geometry from homological algebra is, to me at least, a very attractive one. Certainly, the agreement between spaces of stability conditions and moduli spaces of conformal field theories is impressive enough to suggest that stability conditions do indeed capture some part of the mathematics of string theory. My own feeling is that at some point in the near future the notion of a stability condition will be subsumed into some more satisfactory framework.

The detailed contents of the paper are as follows. Section 2 is about moduli spaces of  superconformal field theories.
This section does not constitute rigorous mathematics but rather provides  background to what follows. Section 3 contains
the basic definitions concerning stability conditions; it is effectively a summary of the contents of  \cite{Br1}.
Section 4 lists the known examples of spaces of stability conditions on smooth projective varieties.
Section 5 contains material on t-structures and tilting; this is used in 
Section 6 to describe other examples of spaces of stability conditions relating to local Calabi-Yau varieties and certain non-commutative algebras. The final, rather speculative section is concerned with extra structures that may exist on spaces of stability conditions.


\section{Background from string theory}

The main players in the mirror symmetry story are N=2 superconformal field theories, henceforth simply referred to as SCFTs. We mathematicians have little intuition for what these are, and until we do mirror symmetry will remain a mystery. This section contains a collection of statements about SCFTs that I have learnt more or less by rote by reading physics papers and talking to physicists.  I do not claim to have any idea of what an SCFT really is and all statements in this section should be taken with a hefty pinch of salt.  Nonetheless, I believe that in the not-too-distant future mathematicians will be able to define SCFTs and prove (or disprove) such statements, and I hope that in the meantime the general picture will be of some value in orienting the reader for what follows. Bourbakistes should skip to Section 3.

\subsection{Abstract SCFTs}

An SCFT has associated to it a pair of topological conformal field theories (TCFTs), the so-called A and B models. It is now fairly well-understood mathematically that a TCFT correpsonds  to a Calabi-Yau $A_{\infty}$ category. The branes of the theory are the
objects of this category. For more precise statements the reader can consult Costello's paper \cite{Co}. However it is possible that TCFTs that arise as topological twistings of SCFTs have some additional properties. Recall that any $A_{\infty}$ category has a minimal model obtained by taking cohomology with respect to the differential $d=m_1$. One thing we shall always require is that this underlying cohomology category is triangulated; as explained in \cite{BK} this is a kind of completeness condition.

The moduli space $\M$ of SCFTs up to isomorphism is expected to be a reasonably well-behaved complex space. Studying this space $\M$ is a good way to try to understand what sort of information an SCFT contains. There is an even better behaved space $\N$ lying over $\M$ which physicists often refer to as Teichm{\"u}ller space, consisting of SCFTs  with some sort of framing. This is the place for example where ratios of central charges of branes are well-defined (see Section \ref{pi} below). One might perhaps hope that $\N$ is a complex manifold with a discrete group action whose quotient is $\M$.

There is an involution on the space $\M$, called the mirror map, induced by an involution of the $N=2$ superconformal algebra. This involution has the effect of exchanging the A and B models. The mirror map permutes the connected components of $\M$ in a non-trivial way.
At general points of $\M$ there are two foliations  corresponding to theories with fixed $A$ or fixed $B$ model respectively. Clearly these must be exchanged by the mirror map. In many examples it seems that the Teichm{\"u}ller space $\N$ has a local product structure inducing the two foliations on $\M$. Note that this implies that the two topological twists together determine the underlying SCFT. The suggestion that this last statement might hold in general seems to be extremely unpalatable to most physicists; it would be interesting to see some examples where it breaks down.

\subsection{Sigma models}

Suppose  $X$ is a Calabi-Yau manifold, and for simplicity assume that $X$ is simply-connected and has complex dimension three.
The non-linear sigma model on $X$ defines an SCFT  depending on a complex structure $I$ on $X$ together with a complexified K{\"a}hler class $\beta+i\omega\in H^2(X,\C)$. In fact, while it is expected that the integrals defining the sigma model will converge providing that the K{\"a}hler class $\omega$ is sufficiently positive, this has  not been proved even to physicists' satisfaction.
Ignoring this problem, the sigma model construction defines an open subset of the moduli space $\M$. More precisely a point in $\M$ is defined by data $(X,I,\beta+i\omega)$ considered up to some discrete group action. For example changing the $B$-field $\beta$ by an integral cohomology class does not effect the isomorphism class of the associated SCFT.  Note also that many points of $\M$ will not be defined by any sigma model, and that points in a single connected component of $\M$ may be defined by sigma models on topologically distinct Calabi-Yau  manifolds.

The A and B twists of the sigma model define TCFTs associated to the Calabi-Yau manifold $X$. The corresponding $A_{\infty}$ categories are  expected to be the derived Fukaya category of the symplectic manifold $(X,\omega)$, suitably twisted by the B-field $\beta$, and an enhanced version of the derived category of coherent sheaves on the complex manifold $(X,I)$ respectively. By enhanced we mean that one should take an $A_\infty$ category whose underlying cohomology category is the usual derived category of coherent sheaves.
If two sigma models corresponding to geometric data $(X_j,I_j,\beta_j+i\omega_j)$ are exchanged by the mirror map one therefore has equivalences
\begin{eqnarray*}
\D^b\Coh (X_1,I_1) &\isom &\D^b\Fuk (X_2,\beta_2+i\omega_2) \\
 \D^b\Fuk (X_1,\beta_1+i\omega_1)&\isom &\D^b\Coh (X_2,I_2).\end{eqnarray*}
This is Kontsevich's homological mirror symmetry proposal \cite{Kon}.

In the neighbourhood of a sigma model defined by data $(X,I,\beta+i\omega)$ the two foliations of $\N$ correspond to varying the complex structure $I$ or the complexified K{\"a}hler class $\beta+i\omega$. In the simply-connected threefold case it is expected that the corresponding leaves in $\M$ are the spaces $\M_{\C}(X)$ of complex structures on $X$ up to diffeomorphism, and the so-called stringy K{\"a}hler moduli space $\M_K(X)$. This latter space has no mathematical definition at present; this is one of the motivations for introducing stability conditions.
The mirror map induces isomorphisms
\[\M_{\C}(X_1)\isom\M_{K}(X_2), \qquad \M_{K}(X_1) \isom\M_{\C}(X_2),\]
and these are often taken as the definition of the K{\"a}hler moduli spaces.

\subsection{K3 surface case}
\label{k3}
The component of the moduli space $\M$ containing sigma models on a K3 surface $X$ was described explicitly by Aspinwall and Morrison \cite{AM}.
It thus gives a good example to focus on, although in some ways it is confusingly different to the picture expected for
a simply-connected  threefold.

The total integral cohomology $\Lambda=H^*(X,\Z)$ has a natural integral symmetric form $(-,-)$
first introduced by Mukai defined by
\[\big( (r_1,D_1,s_1),(r_2,D_2,s_2)\big)=D_1\cdot D_2-r_1 s_2 -r_2 s_1.\]
The resulting lattice $\Lambda$ is even and non-degenerate
and has  signature $(4,20)$. Define
\[\Delta(\Lambda)=\{\delta\in\Lambda:(\delta,\delta)=-2\}.\]
The group $O(\Lambda)$ of
isometries of $\Lambda$ has an index two subgroup $O^+(\Lambda)$
consisting of isometries which preserve the orientation of positive
definite four-planes in $\Lambda\tensor\R$.

The sigma models on $X$ all lie in one connected component $\M_{K3}$ of the moduli space $\M$. The corresponding connected component $\N_{K3}$ of the
Teichm{\"u}ller space $\N$ referred to above is the set of pairs
\[(\Omega,\mho)\in \PP(\Lambda\tensor\C)\]
satisfying the relations
\begin{eqnarray*}
&(\Omega,\Omega)=(\mho,\mho)=(\Omega,\mho)=(\Omega,\bar{\mho})=0, \\
&(\Omega,\bar{\Omega})>0,\quad(\mho,\bar{\mho})>0,
\end{eqnarray*}
with the extra condition
\[\tag{$*$}\text{there is no }\delta\in\Delta(\Lambda)\text{ such that }(\Omega,\delta)=0=(\mho,\delta).\]
The moduli space $\M_{K3}$ is the discrete group quotient $\N_{K3}/O^+(\Lambda)$. Note that $\N_{K3}$ does indeed have a local product structure in this case.
 The mirror map preserves the connected component $\M_{K3}$ (this is the statement that K3 surfaces are self-mirror) and simply exchanges
 $\Omega$ and $\mho$.

In fact Aspinwall and Morrison do not impose the last condition $(*)$. 
In the physics there is indeed some sort of theory existing at the points where $(*)$ fails, but the brane with the corresponding charge $\delta$ has become massless, so that the strict SCFT description breaks down, and non-perturbative corrections have to be taken into account. In any case, it is clear that for our purposes it is important to leave out the hyperplanes where $(*)$ fails, since one obtains interesting transformations by taking monodromy around them.

The sigma model defined by a complex structure $I$ and a complexified K{\"a}hler class $\beta+i\omega$ on $X$
corresponds to the pair $(\Omega,\mho)$ where $\Omega$ is the class of a holomorphic two-form on $(X,I)$ and
\[\mho=[e^{\beta+i\omega}]=[(1,\beta+i\omega,\frac{1}{2}(\beta+i\omega)^2)]\in\PP (\Lambda\tensor\C).\]
It is important to note that in the K3 case, in contrast to the case of a simply-connected threefold, the set of choices of this data does not define an open subset of the moduli space $\M_{K3}$.  In fact the point of $\M_{K3}$ defined by a pair $(\Omega,\mho)$ can only come from a sigma model of the usual sort  if there is a hyperbolic plane in $\Lambda$ orthogonal to $\Omega$. To obtain an open set in $\M_{K3}$ one must also include sigma models defined using  Hitchin's generalised complex structures \cite{Hi,KL}. 
As explained by Huybrechts \cite{Hu} one then obtains the entire space $\M_{K3}$ as the moduli space of
generalised K3 structures on $X$.

\subsection{Stability for $D$-branes}
\label{pi}

Suppose that one is given an SCFT together with one of its topological twists. As above, this twist corresponds mathematically to a Calabi-Yau $A_{\infty}$ category $\D$. Some deformations of the SCFT will induce deformations of $\D$, but there is a leaf $\LL$ in the Teichm{\"u}ller space $\N$ whose points all have the same topological twist. The notion of a stability condition comes about by asking what the significance of these extra parameters is for the category $\D$.

Douglas argued \cite{Do2,Do3} that at each point on the leaf $\LL\subset\N$ there is a full subcategory of semistable objects or BPS branes in the category $\D$. On a $\C-$bundle over this leaf one can also assign complex numbers $Z(E)$ called central charges to all objects $E\in\D$. Moreover, for a BPS brane $E$ the central charge $Z(E)$ is always nonzero, and there is a distinguished choice of phase
\[\phi(E)=\frac{1}{\pi}\arg Z(E)\in\R.\]
Axiomatising the properties of these subcategories $\P(\phi)\subset\D$ of BPS branes of phase $\phi$ leads to the definition of a stability condition given in the next section.

A good example to bear in mind for heuristic purposes is the case when the SCFT is a sigma model on a Calabi-Yau threefold $(X,I,\beta+i\omega)$ and \[\D=\D^b\Fuk(X,\beta+i\omega)\]  is the A-model topological twist. Then the leaf $\LL\subset\N$ is expected to be the space of complex structures on $X$ up to diffeomorphisms isotopic to the identity. On a $\C^*$-bundle over $\LL$ there is a well-defined holomorphic three-form $\Omega$.
Given this choice of $\Omega$, every Lagrangian submanifold $L\subset X$ has an associated complex number
\[Z(L)=\int_L\Omega\in\C.\]
A Lagrangian submanifold $L\subset X$ is said to be special of phase $\phi\in\R/2\Z$ if
\[\Omega|_L=\exp(i\pi\phi) \vol_L,\]
where $\vol_L$ is the volume form on $L$. Note that one then has
\[Z(L)=\int_L e^{i\pi\phi} \vol_L \in\R_{>0}\exp({i\pi\phi}).\]

To be slightly more precise, the objects of the Fukaya category are graded Lagrangians. This means that the phases $\phi$ of the special Lagrangians can be lifted to elements of $\R$ up to an overall integer indeterminacy which can be eliminated by passing to a $\Z-$cover. Thus corresponding to each point of a $\C-$bundle over $\LL$ and each $\phi\in\R$ there is a full subcategory $\P(\phi)\subset\D$ of special Lagrangians of phase $\phi$.

If one varies the complex structure on $X$ and hence the holomorphic three-form $\Omega$, the numbers $Z(L)$ and the subcategories of special Lagrangians will change. This process was studied by Joyce \cite{Jo} and the analogy with variations of stability in algebraic geometry was explained by Thomas \cite{Th}. These ideas provided the basis for Douglas' work on $\Pi$-stability for D-branes in a general SCFT context.


\section{Basic definitions}
This section is a summary of the essential definitions and results from \cite{Br1}. Proofs and further details can be found there.

\subsection{Stability conditions}
For the rest of the paper  $\D$ will denote a triangulated category and $K(\D)$ its Grothendieck group. In the previous section we worked with an $A_\infty$ category whose cohomology category was triangulated; no such enhancement will be necessary for what follows. 

\begin{defn}
\label{pemb}
A stability condition $\sigma=(Z,\P)$ on  $\D$
consists of
a group homomorphism
$Z\colon K(\D)\to\C$ called the central charge,
and full additive
subcategories $\P(\phi)\subset\D$ for each $\phi\in\R$,
satisfying the following axioms:
\begin{itemize}
\item[(a)] if $0\neq E\in \P(\phi)$ then $Z(E)=m(E)\exp(i\pi\phi)$ for some
 $m(E)\in\R_{>0}$,
\item[(b)] for all $\phi\in\R$, $\P(\phi+1)=\P(\phi)[1]$,
\item[(c)] if $\phi_1>\phi_2$ and $A_j\in\P(\phi_j)$ then $\Hom_{\D}(A_1,A_2)=0$,
\item[(d)] for each nonzero object $E\in\D$ there is a finite sequence of real
numbers
\[\phi_1>\phi_2> \cdots >\phi_n\]
and a collection of triangles
\[
\xymatrix@C=.5em{
0_{\ } \ar@{=}[r] & E_0 \ar[rrrr] &&&& E_1 \ar[rrrr] \ar[dll] &&&& E_2
\ar[rr] \ar[dll] && \ldots \ar[rr] && E_{n-1}
\ar[rrrr] &&&& E_n \ar[dll] \ar@{=}[r] &  E_{\ } \\
&&& A_1 \ar@{-->}[ull] &&&& A_2 \ar@{-->}[ull] &&&&&&&& A_n \ar@{-->}[ull] 
}
\]
with $A_j\in\P(\phi_j)$ for all $j$.
\end{itemize}
\end{defn}

\begin{remark}
The central charge part of the definition of a stability condition is mathematically quite bizarre.
For example it means that even if one takes the derived category of a variety defined over a field of positive characteristic one still
obtains a complex manifold as the space of stability conditions.
Omitting the central charge in Definition \ref{pemb} gives the more natural notion of a slicing. Note however that the  group of orientation-preserving homeomorphisms of the circle acts on the set of slicings of $\D$ by relabelling the phases $\phi$, so that spaces of slicings will always be hopelessly infinite-dimensional.
\end{remark}

Given a stability condition $\sigma=(Z,\P)$ as in the definition, each subcategory
$\P(\phi)$ is abelian. The nonzero objects of $\P(\phi)$ are said to be
{semistable of phase $\phi$ in $\sigma$}, and the simple objects
of $\P(\phi)$ are said to be stable.
It follows from the other axioms that the decomposition of
an
object $0\neq E\in\D$ given by
axiom (d) is uniquely defined up to isomorphism.
Write
$\phi^+_{\sigma}(E)=\phi_1$ and
$\phi^-_{\sigma}(E)=\phi_n$. 
The {mass} of $E$ is defined to be
the positive real number
$m_{\sigma}(E)=\sum_i |Z(A_i)|$.

For any interval $I\subset\R$, define $\P(I)$ to be the extension-closed
subcategory of $\D$ generated by the subcategories $\P(\phi)$ for
$\phi\in I$. Thus, for example, the full subcategory $\P((a,b))$
consists of the zero objects of $\D$ together with those
objects $0\neq E\in\D$ which satisfy
$a<\phi_\sigma^-(E)\leq\phi_\sigma^+(E)<b$.

In order for Theorem \ref{lasty} below to hold it is necessary to restrict attention to stability conditions
satisfying an extra technical condition. A stability condition is called locally finite if
there is some $\epsilon>0$ such that each quasi-abelian category
$\P((\phi-\epsilon,\phi+\epsilon))$ is of finite length. For details see \cite{Br1}.

\subsection{Spaces of stability conditions}
\label{bird}
We always assume that our triangulated categories are essentially small, i.e. equivalent to a category with a set of objects.
Write $\Stab(\D)$ for the set of locally-finite stability conditions on a fixed
triangulated category $\D$. It has a natural topology induced by the
metric
\[d(\sigma_1,\sigma_2)=\sup_{0\neq E\in\D}
\bigg\{|\phi^-_{\sigma_2}(E)-\phi^-_{\sigma_1}(E)|,|\phi^+_{\sigma_2}(E)
-\phi^+_{\sigma_1}(E)|
,|\log \frac{m_{\sigma_2}(E)}{m_{\sigma_1}(E)}|\bigg\}
\in[0,\infty].\]

\begin{remark}
\label{closed}
Note that the set of stability conditions $\sigma\in\Stab(\D)$ for which a given object $E\in\D$ is semistable is closed. Indeed, a nonzero object $E\in\D$ is semistable in a stability condition $\sigma$ precisely if $\phi^+_\sigma(E)=\phi_\sigma^-(E)$.
\end{remark}

\begin{remark}
Is it possible to use properties of the above metric to show that spaces of stability conditions are always contractible? Such a result would have non-trivial implications for computing groups of autoequivalences of derived categories. There are no known counterexamples.
\end{remark}

It is clear that there is a forgetful map
\[\ZZ\colon\Stab(\D)\lra \Hom_{\Z}(K(\D),\C)\]
sending a stability condition to its central charge.
The following  result was proved in \cite{Br1}. Its slogan is that deformations of the central charge
lift uniquely to deformations of the stability condition. One has to be a little bit careful because in general (in fact even for $\D=\D^b\Coh(X)$ with $X$ an elliptic curve) the group $K(\D)$ has infinite rank.

\begin{thm}
\label{lasty}
Let $\D$ be a triangulated category.
For each connected component $\Stab^*(\D)\subset\Stab(\D)$
there is a linear subspace $V\subset \Hom_{\Z}(K(\D),\Z)$,
with a well-defined linear topology, such that the restricted map
$\ZZ\colon\Stab^*(\D)\to \Hom_{\Z}(K(\D),\C)$
is a local homeomorphism onto an open subset of $V$.
\end{thm}

In the cases of most interest for us we can get round the problem of $K(\D)$ being of infinite rank as follows.

\begin{defn}
If $X$ is a smooth complex projective manifold we write $\Stab(X)$ for the set of locally-finite stability conditions on $\D=\D^b\Coh(X)$ for which
the central charge $Z$ factors via the Chern character map
\[
\xymatrix@C=1.5em{ K(\D)\ar[d]^{\ch}\ar[rr]^{Z} && \C\\
H^*(X,\QQ)\ar@{-->}[urr] 
}
\]
In physics language we are insisting that the central charge of a brane depends only on its topological charge.
\end{defn}

Theorem \ref{lasty} immediately implies that $\Stab(X)$ is a finite-dimensional complex manifold.

\begin{remark}
\label{pc}
There are other more natural ways to get round the problem of $K(\D)$ being infinite rank which will almost certainly play a r\^ole in the further development of the theory. One possibility would be to take some homology theory $H(\D)$ associated to $\D$ for which there is a natural Chern character map $K(\D)\to H(\D)$, and then insist that the central charge $Z$ factors via $H(\D)$.  For example one could take periodic cyclic homology. In the case when $\D=\D^b\Coh(X)$, with $X$  a smooth projective variety, it is known \cite{We} that this homology theory coincides as a vector space with the de Rham cohomology of $X$, and in particular is finite-dimensional. 
But in the absence of any real understanding of why this would be particularly a sensible thing to do we shall stick with the ad hoc definition above.
\end{remark}

\subsection{Group actions}
The space of stability conditions of any triangulated category has some natural group actions which will be important in what follows.

\begin{lemma}
\label{groupactions}
The space $\Stab(\D)$ carries a right action of the
group
$\grp$, the universal
covering space of $\GL^+(2,\R)$, and a left action of the 
group
$\Aut(\D)$ of exact autoequivalences of $\D$. These two actions
commute.
\end{lemma}

\begin{pf}
First note that the group $\grp$ can be thought of as the
set of pairs $(T,f)$ where $f\colon \R\to\R$ is an
increasing map with $f(\phi+1)=f(\phi)+1$, and  $T\colon \R^2\to\R^2$
is an orientation-preserving linear isomorphism, such that the induced maps on
$S^1=\R/2\Z=(\R^2\setminus\{0\})/\R_{>0}$ are the same.

Given a stability condition $\sigma=(Z,\P)\in\Stab(\D)$, and a
pair $(T,f)\in\grp$, define a new stability
condition $\sigma'=(Z',\P')$ by setting $Z'=T^{-1}\circ Z$ and
$\P'(\phi)=\P(f(\phi))$. Note that the semistable objects of the
stability conditions $\sigma$
and $\sigma'$ are the same, but the phases have been
relabelled.

For the second action, note that an element
$\Phi\in\Aut(\D)$ induces an automorphism $\phi$ of $K(\D)$. If
$\sigma=(Z,\P)$ is a stability condition on $\D$ define $\Phi(\sigma)$
to be the stability condition $(Z\circ\phi^{-1},\P')$, where
$\P'(t)$=$\Phi(\P(t))$. 
\end{pf}

Neither of the two group actions of Lemma \ref{groupactions} will be free in general. In particular, if $\sigma=(Z,\P)$ is a stability
condition in which the image of the central charge $Z\colon K(\D)\to\C$ lies on a real line in $\C$
then $\sigma$ will be fixed by some subgroup of $\grp$. However there is a subgroup $\C\subset\grp$ which does act freely.
If $\lambda\in\C$ then $\lambda$ sends a stability condition $\sigma=(Z,\P)$ to the stability condition $\lambda(\sigma)=(Z',\P')$ where
$Z'(E)=e^{-i\pi\lambda}Z(E)$ and $\P'(\phi)=\P(\phi+\Re(\lambda))$.
Note that for any integer $n$ the action of the shift functor $[n]$ on $\Stab(\D)$ coincides with the action of $n\in\C$.

\remark
\label{birdtwo}
Return for a moment to the discussion of Section \ref{pi} in which \[\D=\D^b\Fuk(X,\beta+i\omega).\]  The action of $\C$ on $\Stab(\D)$ clearly corresponds to rotating the holomorphic three-form $\Omega$. It also seems reasonable to guess that the action of $\Aut(\D)$ on $\Stab(\D)$ corresponds to the discrete group quotient $\N\to\M$.
Thus we might expect an embedding of the complex moduli space $\M_{\C}(X)$ in the double quotient
\[\Aut(\D)\backslash\Stab(\D)/\C.\]
The mirror statement is that if  $X$ is a Calabi-Yau with a given complex structure and $\D=\D^b\Coh(X)$ then the above quotient contains
the stringy K{\"a}hler moduli space $\M_K(X)$. In the next section we will examine this suggestion in some simple examples.


\section{Compact examples}

In this section we review the known examples of stability conditions on smooth projective varieties. The only Calabi-Yau examples are elliptic curves and K3 and abelian surfaces. 

\subsection{Elliptic curves}
Let $X$ be a complex projective curve of genus one. It was shown in \cite{Br1} that the action of $\grp$ on $\Stab(X)$
is free and transitive. Thus
\[\Stab(X)\isom\grp\isom\C\times\HH\]
where $\HH\subset\C$ is the upper half-plane. Quotienting by the group of autoequivalences of $\D=\D^b\Coh(X)$ gives
\[\frac{\Stab(X)}{\Aut(\D)}\isom \frac{\GL^+(2,\R)}{\PSL(2,\Z)}\]
which is thus a $\C^*-$bundle over the modular curve $\HH/\PSL(2,\Z)$. In fact this is the $\C^*-$bundle parameterising
equivalence classes of data consisting of a complex structure on $X$ together with a non-zero holomorphic $1$-form.

According to Remark \ref{birdtwo}  we expect an inclusion of the stringy K{\"a}hler moduli space in the quotient
\[\Aut(\D)\backslash\Stab(X)/\C\isom \HH/\PSL(2,\Z).\]
We also know that
\[\M_K(X)=\M_{\C}(\check{X}),\]
and the since tori are self-mirror the latter is just the moduli of complex structures on
a two torus.
Thus we obtain perfect agreement in this case.

\begin{remark}
The calculation of the space of stability conditions on an elliptic curve has been generalised by Burban and Kreussler to include irreducible singular curves of arithmetic genus one \cite{Burb}. The resulting space of stability conditions and the quotient by the group of autoequivalences is the same as in the smooth case. 
\end{remark}

\begin{remark}
It is possible that spaces of stability conditions can be (partially) compactified by adding non-locally-finite stability conditions. In the example of the elliptic curve the non-locally-finite stability conditions up to the action of $\C$ are parameterised by $\R\setminus\QQ$. 
It might be interesting to think this point through in some other examples.
\end{remark}

\subsection{K3 surfaces}
Let $X$ be an algebraic K3 surface and set $\D=\D^b\Coh(X)$. I use the notation introduced in Section \ref{k3}; in particular $\Lambda$ denotes the lattice $H^*(X,\Z)$ equipped with the Mukai symmetric form. Let $\Omega\in H^2(X,\C)$ be the class of a nonzero holomorphic two-form on $X$; we consider $\Omega$ as an element of $\Lambda\tensor\C$. 
The
sublattice
\[\N(X)=\Lambda\cap
\Omega^{\perp}\subset H^*(X,\Z)\] can be identified with $\Z\oplus\NS(X)\oplus\Z$
and has signature $(2,\rho)$, where
$1\leq\rho\leq 20$ is the Picard number of $X$.
Write $O^+(\Lambda,\Omega)$ for the subgroup of $O^+(\Lambda)$ consisting of isometries which preserve the class $[\Omega]\in\PP(\Lambda\tensor\C)$. Any such isometry restricts to give an isometry of $\N(X)$.
Set $\Delta(\Lambda,\Omega)=\Delta(\Lambda)\cap\Omega^{\perp}$ and for each
$\delta\in\Delta(\Lambda,\Omega)$ let
\[\delta^{\perp}=\{\mho\in\N(X)\tensor\C:(\mho,\delta)=0\}\subset
\N(X)\tensor\C\]
be the
corresponding complex hyperplane.

Define an open subset $\Pol(X)\subset\N(X)\tensor\C$
consisting of vectors $\mho\in\N(X)\tensor\C$ whose real and imaginary parts span a positive definite two-plane in $\N(X)\tensor\R$.
Taking orthogonal bases in the two-plane shows that $\Pol(X)$ is a $\GL(2,\R)$-bundle over the set
\[\Q(X)=\{\mho\in \PP(\N(X)\tensor\C):(\mho,\mho)=0\text{ and }(\mho,\bar\mho)>0\}.\]
These spaces $\P(X)$ and $\Q(X)$ have two connected components that are exchanged
by complex conjugation. Let $\Pol^+(X)$ and $\Q^+(X)$ denote the ones containing vectors of the form $(1,i\omega,-\omega^2/2)$ with $\omega\in\NS(X)$
the class of an ample line bundle.

The Mukai vector of an object $E\in\D(X)$ is defined to be
\[v(E)=\ch(E)\sqrt{\td(X)}\in \N(X),\]
where $\ch(E)$ is the Chern character of $E$ and $\td(X)$ is the Todd class of $X$.
The fact that the Mukai form is non-degenerate means that for any $\sigma=(Z,\P)\in\Stab(X)$ we can write the central charge $Z$ in the form
\[Z(E)=(\pi(\sigma),v(E))\]
for some vector $\mho=\pi(\sigma)\in\N(X)$. This defines a map $\pi\colon\Stab(X)\to\N(X)\tensor\C$.

It was proved in \cite{Br2}
that there is a connected component $\Stab^\dagger(X)\subset\Stab(X)$
that is mapped by $\pi$ onto the open subset
\[\Pol_0^+(X)=\Pol^+(X)\setminus\bigcup_{\delta\in\Delta(\Lambda,\Omega)}\delta^{\perp}
\subset\N(X)\tensor\C.\]
Moreover, the induced map
$\Stab^\dagger(X)\to\Pol^+_0(X)$
is a covering map, and if $\Aut^\dagger(\D)$ is the subgroup of $\Aut(\D)$ preserving the connected component $\Stab^\dagger(X)$,
then \[\frac{\Stab^\dagger(X)}{\Aut^\dagger(\D)}\isom \frac{\Pol^+_0(X)}{O^+(\Lambda,\Omega)}.\]
Comparing this result with the discussion in Section \ref{k3} we see that the space
\[\Aut^\dagger(\D)\backslash\Stab^\dagger(X)/\C\]
agrees closely with the leaf $\M_K(X)\subset\M_{K3}$ corresponding
to the fixed holomorphic two-form $\Omega$, but that there are two differences.

Firstly, the space $\M_{K}(X)$ consists of points $\mho$ in the projective space of
\[\Omega^{\perp}\cap{\bar{\Omega}}^{\perp}=\bigoplus_{p=0}^2 H^{p,p}(X)\subset H^*(X,\C),\]
whereas the stability conditions space
only sees the algebraic part $\N(X)\tensor\C$. This problem could presumably be fixed by changing the definition of the map $Z$ as suggested in Remark \ref{pc}.

The second point is that the vectors $\mho$ in $\M_K(X)$ satisfy $(\mho,\mho)=0$ whereas the space of stability conditions has no such normalisation. This means that $\Stab(X)$ is one complex dimension larger than one would otherwise expect. This may seem a minor point but is in fact very important, being the first glimpse of Hodge theoretic restrictions on the central charge. We discuss such restrictions further in Section \ref{font}.

\subsection{Other examples}

The only other compact varieties for which stability conditions are known are abelian surfaces and varieties of dimension one. The case of abelian surfaces is very similar to the K3 case but easier; it is covered in \cite{Br2}.
It might also be possible to calculate the space of stability conditions on  higher dimensional abelian varieties,
 although I have not thought about this in any detail. 
Turning to non-Calabi-Yau examples S. Okada \cite{Ok} proved that
\[\Stab(\PP^1)\isom \C^2\]
and E. Macri \cite{Ma} proved that for any curve $X$ of genus $g\geq 2$ one has
\[\Stab(X)\isom \grp\isom\C\times\HH.\]
Of course it is also possible that the definition of stability condition needs to be changed for non Calabi-Yau categories to take account of the non-trivial Serre functor in some way.

It would be extremely interesting to calculate the space of stability conditions on a compact Calabi-Yau threefold such as the quintic. Unfortunately we don't know enough about coherent sheaves on threefolds to be able to do this. In particular it would be useful to know the set of Chern characters of Gieseker stable bundles. In fact it is not known how to write down a single stability condition on a Calabi-Yau threefold.
Nonetheless, physicists have used mirror symmetry together with a certain amount of guesswork to make some nontrivial computations \cite{AD,AK}.

One way to proceed for the quintic would be to construct its derived category via matrix factorisations as in \cite{Or} and then to somehow construct the stability condition corresponding to the Gepner point in the stringy K{\"a}hler moduli space. This seems to be an interesting project. For an example of a stability condition on a category of matrix factorisations see \cite{Sai,Ta}.


\section{T-structures and tilting}
\label{rowan}

In this section I explain the connection between stability conditions and t-structures. This is the way
stability conditions are constructed in practice. I also explain how the method of tilting can be
used to give a combinatorial description of certain spaces of stability conditions. This technique will be
applied in the next section to describe  spaces of stability conditions on some non-compact Calabi-Yau varieties.

\subsection{Stability conditions and t-structures}

A  stability function on
an abelian  category $\A$ is defined to be  a group homomorphism $Z\colon K(\A)\to\C$
such that
\[0\neq E\in\A \implies Z(E)\in\R_{>0}\,\exp({i\pi\phi(E)})\text{ with
}0<\phi(E)\leq 1.\]
The real number $\phi(E)\in(0,1]$ is called the phase of the object $E$.

A nonzero object $E\in\A$ is said to be
semistable
with respect to $Z$
if every subobject $0\neq A\subset E$ satisfies $\phi(A)\leq\phi(E)$.
The stability function $Z$ is said to have the Harder-Narasimhan property if
every nonzero object $E\in\A$ has
a finite filtration
\[0=E_0\subset E_1\subset \cdots\subset E_{n-1}\subset E_n=E\]
whose factors $F_j=E_j/E_{j-1}$ are semistable objects of $\A$ with
\[\phi(F_1)>\phi(F_2)>\cdots>\phi(F_n).\]

Given a stability condition $\sigma=(Z,\P)$ on a triangulated category $\D$ the full subcategory $\A=\P((0,1])\subset\D$ is the heart of a bounded t-structure on $\D$.  It follows that $\A$ is an abelian category and we can identify its Grothendieck group $K(\A)$ with $K(\D)$. We call $\A$ the heart of the stability condition $\sigma$. The central charge $Z$ defines a stability function on
$\A$, and the decompositions of axiom (d) give Harder-Narasimhan filtrations.

Conversely, given a bounded t-structure on $\D$ together with a stability function $Z$ on its heart $\A\subset\D$, we can define subcategories
$\P(\phi)\subset\A\subset\D$ to be the semistable objects in $\A$ of phase $\phi$ for each $\phi\in(0,1]$. Axiom (b) then fixes $\P(\phi)$ for all $\phi\in\R$, and identifying $K(\A)$ with $K(\D)$ as before, we obtain a stability condition $\sigma=(Z,\P)$ on $\D$.
Thus we have

\begin{prop}
\label{pg}
To give a stability condition on a triangulated category $\D$ is
equivalent to giving a bounded t-structure on $\D$ together with a 
stability function on its heart with the Harder-Narasimhan property.
\end{prop}

Let $\D$ be a triangulated category and suppose $\A\subset\D$ is the heart of a bounded t-structure on $\D$. We write $U(\A)$ for the subset of $\Stab(\D)$ consisting of stability conditions with heart $\A$. In general this subset could be empty. Suppose though that $\A$ is a finite length
 category with finitely many isomorphism classes of simple objects $S_1,\cdots, S_n$. The Grothendieck group $K(\D)=K(\A)$ is then a free abelian group on the generators $[S_i]$. 
Set \[H=\{r\exp(i\pi\phi):r\in\R_{>0}\text{ and }0<\phi\leq 1\}\subset\C.\]
According to Proposition \ref{pg}
we can define a stability condition
on $\D$ with heart $\A$ by choosing a central charge $Z(S_i)\in H$ for each $i$ and extending linearly to give a map
$Z\colon K(\D)\to\C$. This argument gives

\begin{lemma}
\label{ram}
Let $\A\subset\D$ be the heart of a bounded t-structure on $\D$ and suppose $\A$ is finite length with $n$ simple objects $S_1,\cdots,S_n$. Then the
subset $U(\A)\subset \Stab(\D)$ consisting of stability conditions with heart $\A$ is isomorphic to $H^n$.
\end{lemma}

The next step is to understand stability conditions on the boundary of the region $U(\A)$ described above.
To do this we need the method of tilting.

\subsection{Tilting}

In the level of generality we shall need, tilting was introduced
by Happel, Reiten and Smal{\o} \cite{HRS}, although the name and the basic idea go back to a paper of Brenner and Butler \cite{BB}.

\begin{defn}
\label{tors} A torsion pair in an abelian category $\A$ is a pair of full
subcategories $(\T,\F)$ of $\A$ which satisfy $\Hom_{\A}(T,F)=0$ for $T\in
\T$ and $F\in\F$, and such that every object $E\in\A$ fits into a  short
exact sequence
\[0\lra T\lra E\lra F\lra 0\] for some pair of objects $T\in\T$ and
$F\in  \F$.
\end{defn}

The objects of $\T$ and $\F$ are called torsion and torsion-free. The  following result \cite[Proposition 2.1]{HRS} is
easy to check.

\begin{prop}(Happel, Reiten, Smal{\o})
Suppose $\A\subset \D$ is the heart of a bounded t-structure on a triangulated category
$\D$. Given an object $E\in\D$ let $H^i(E)\in\A$ denote the $i$th cohomology
object of $E$ with respect to this t-structure. Suppose $(\T,\F)$ is a
torsion pair in $\A$. Then the full subcategory
\[\A^{\sharp}=\big\{E\in \D:H^i(E)=0\text{ for }i\notin\{-1,0\},
H^{-1}(E)\in\F\text{ and } H^{0}(E)\in\T\big\}\] is the heart of a bounded
t-structure on $\D$.\qed
\end{prop}

Recall that a bounded t-structure on $\D$ determines and is determined by its heart $\A\subset\D$.
In the situation of the Lemma one says that the the subcategory $\A^{\sharp}$ is
obtained from the subcategory $\A$ by tilting with respect to the
torsion pair $(\T,\F)$. In
fact one could equally well consider $\A^{\sharp}[-1]$ to be the tilted
subcategory; we shall be more precise about this where necessary. Note that the
pair $(\F[1],\T)$ is a torsion pair in $\A^{\sharp}$ and that tilting with
respect to this pair gives back the original subcategory $\A$ with a shift.

Suppose $\A\subset \D$ is the heart of a bounded t-structure and is a
finite length abelian category.  Given a simple object $S\in \A$ define
$\langle S \rangle\subset \A$ to be the full subcategory consisting of objects $E\in\A$ all
of whose simple factors are isomorphic to $S$. One can either view $\langle S \rangle$ as
the torsion part of a torsion theory on $\A$, in which case the torsion-free
part is
\[\F=\{E\in \A:\Hom_{\A}(S,E)=0\},\]
or as the torsion-free part, in which case the torsion part is
\[\T=\{E\in\A:\Hom_{\A}(E,S)=0\}.\]
The corresponding tilted subcategories are defined to be
\begin{eqnarray*}
\L_S \A &=& \{E\in\D:H^i(E)=0\text{ for
}i\notin\{0,1\},H^{0}(E)\in\F\text{ and }H^1(E)\in\langle S \rangle\}  \\
\RR_S \A &=& \{E\in \D:H^i(E)=0\text{ for
}i\notin\{-1,0\},H^{-1}(E)\in\langle S \rangle\text{ and }H^0(E)\in\T\}.
\end{eqnarray*}

We can now return to stability conditions. Suppose we are in the situation of Lemma \ref{ram} and $\sigma=(Z,\P)$ is a stability condition in the boundary of the region $U(\A)$. Then there is some $i$ such that $Z(S_i)$ lies on the  real axis. Forgetting about higher codimension phenomena for now let us assume that $\Im Z(S_j)>0$ for every $j\neq i$. Since each object $S_i$ is stable for all stability conditions in $U(\A)$, by Remark \ref{closed}, each $S_i$ is at least semistable in $\sigma$, and hence $Z(S_i)$ is nonzero.
The following result is easily checked.

\begin{lemma}
\label{tilting}
In the situation of Lemma \ref{ram} suppose $\sigma=(Z,\P)\in\Stab(\D)$ lies on a unique codimension one boundary of the region $U(\A)$.
Then either $Z(S_i)\in\R_{<0}$ for some $i$ and a neighbourhood of $\sigma$ is contained in $U(\A)\cup U(\L_{S_i}(\A))$,
or $Z(S_i)\in\R_{>0}$ for some $i$ and a neighbourhood of $\sigma$ is contained in $U(\A)\cup U(\RR_{S_i}(\A))$.
\end{lemma}

In general the tilted subcategories $\L_{S_i}(\A)$ and $\RR_{S_i}(\A)$ need not be finite length and so we cannot necessarily repeat this process. But in many examples we can. Then we obtain a subset  of $\Stab(\D)$ covered by regions isomorphic to $H^n$, each one corresponding to a given heart $\A\subset\D$, and with different regions glued together along boundaries corresponding to pairs of hearts related by tilts at simple objects.
Thus understanding the algebra of the tilting process can lead to a
combinatorial description for certain  spaces of stability conditions. In the next section we shall see some examples of this.


\section{Non-compact examples}

Given the primitive state of knowledge concerning coherent sheaves on projective varieties of dimension at least three it is natural to
study quasi-projective
 varieties instead.
A particularly amenable class of examples consists of varieties for which there exists a derived equivalence
\[\D^b\Coh(X)\isom\D^b\Mod(B),\]
where $\Mod(B)$ is the category of modules over some non-commutative algebra $B$.
In practice, the non-compact variety $X$ is often the total space of a holomorphic vector bundle on a lower-dimensional variety $Z$; such examples are called local varieties in the physics literature. The derived equivalence is then obtained using the theory of exceptional collections \cite{rud} and the relevant algebras $B$ can be described via a quiver with relations \cite{bon}. 

\subsection{Some local Calabi-Yau varieties}
Suppose $Z$ is a Fano variety and let $X=\omega_Z$ be the total space of the canonical bundle of $Z$ with its projection $\pi\colon X\to Z$. Suppose $Z$ has a full exceptional collection $(E_0,\cdots,E_{n-1})$ such that for all $p\leq 0$ and all $i,j$ one has
\[\Hom^k_Z(E_i,E_j\tensor\omega_Z^p)=0 \text { unless } k= 0.\]
Such collections are called geometric and are known to exist in many interesting examples.
It was shown in \cite{Br3} that there is then an equivalence
\[\D^b\Coh(X)\isom \D^b\Mod(B),\]
as above, where $B$ is the endomorphism algebra
\[B=\End_X\big(\bigoplus_i \pi^* E_i\big).\]
We give some examples.

\begin{example}
Take $Z=\PP^2$ and the geometric exceptional collection \[\OO,\OO(1),\OO(2).\] Then the algebra $B$ is the path algebra of the quiver
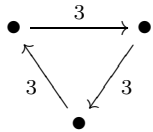
\[
\xymatrix@C=1em{ \bullet\ar[rr]^{3} && \bullet\ar[dl]^{3}\\
&\bullet\ar[ul]^{3} 
}
\]
with commuting relations. The numbers labelling the arrows denote the number of arrows. Commuting relations
 means that if we label the arrows joining a pair of vertices by symbols $x,y,z$ then
the relations are given by $xy-yx, yz-zy$ and $xz-zx$.
\end{example}

\begin{example}
Take $Z=\PP^1\times\PP^1$ and the geometric exceptional collection \[\OO,\OO(1,0),\OO(0,1),\OO(1,1).\]
This leads to a quiver of the form
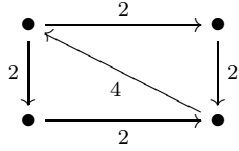
\[
\xymatrix@C=5em{ \bullet\ar[d]_{2}\ar[r]^{2} & \bullet\ar[d]^{2}\\
\bullet\ar[r]_{2} &\bullet \ar[ul]^{4}
}
\]
with some easily computed relations.
Alternatively, one could take the  collection \[\OO,\OO(1,0),\OO(1,1),\OO(2,1),\]
which gives the quiver
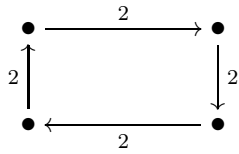
\[
\xymatrix@C=5em{ \bullet\ar[r]^{2} & \bullet\ar[d]^{2}\\
\bullet\ar[u]^{2} &\bullet \ar[l]^{2}
}
\]
\end{example}

It is usually convenient to restrict attention the full subcategory $\D\subset\D^b\Coh(X)$ consisting of objects supported on the zero-section $Z\subset X$.
In a sense this is the interesting part of the category, the rest being rather flabby.
The equivalence above restricts to give an equivalence between $\D$ and the full subcategory of $\D^b\Mod(B)$ consisting of complexes with nilpotent cohomology modules. The standard t-structure on $\D^b\Mod(B)$ then induces a bounded t-structure on $\A\subset \D$ whose heart is equivalent to the finite length category $\Mod_0(B)$ of nilpotent $B$-modules. We call the subcategories $\A\subset \D$ obtained from geometric exceptional collections in this way exceptional.

\subsection{Tilting, mutations and Seiberg duality}

In order to apply the theory of Section \ref{rowan} to the situation described above we need to understand the tilting process for exceptional subcategories.
This problem is closely related to the theory of mutations of exceptional collections, as studied in the Rudakov seminar \cite{rud}. The connection is worked out in detail in \cite{St}, see also  \cite{AMe,BD,He}.
 
In the case when $Z$ satisfies
\[n=\dim K(Z)\tensor\C=1+\dim Z\]
(for example when $Z$ is a projective space) the relationship between tilting and mutations is particularly straightforward and is described precisely in \cite{Br3}.
In particular any tilt of an exceptional subcategory $\A\subset\D$ is the image of an exceptional subcategory by some autoequivalence of $\D$. This means that the tilting process can be continued indefinitely. The combinatorics of the process is controlled by the affine braid group
\[B_n=\big\langle \tau_0,\cdots,\tau_{n-1} \;|\;\tau_i\tau_j\tau_i=\tau_j\tau_i\tau_j\text{ for all }j-i\equiv 1 \mod n\big \rangle.\]
For details on this result we refer the reader to \cite{Br3}; the main input is Bondal and Polishchuk's work \cite{BP}.
In the case $Z=\PP^2$ these results allow one to give a combinatorial description  of a connected component of $\Stab(\D)$.

\begin{thm}[\cite{Br5}]
\label{one}
Set $Z=\PP^2$ and define $\D\subset\D^b\Coh(\omega_Z)$ as above.
Then there is a subset of $\Stab(\D)$ that can be written as a disjoint union of regions
\[\bigsqcup_{g\in G} U(g),\]
where $G=B_3$ is the affine braid group on three strings. 
Each region $U(g)$ is isomorphic to $H^3$ and consists of  stability conditions with a given heart $\A(g)\subset \D$. The closures of
two regions $U(g_1)$ and $U(g_2)$ intersect along
a codimension one boundary  precisely if $g_1 g_2^{-1}=\tau_i^{\pm 1}$ for some $i$.
Each of the categories $\A(g)$ is equivalent to a category of nilpotent representations of a quiver with relations of the form
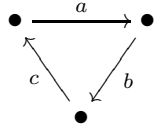
\[
\xymatrix@C=1em{ \bullet\ar[rr]^{a} && \bullet\ar[dl]^{b}\\
&\bullet\ar[ul]^{c} 
}
\]
where the positive integers $a,b,c$ counting the numbers of arrows connecting the vertices always satisfy the Markov equation
\[a^2 + b^2 + c^2 =abc.\]
\end{thm}

The web of categories indexed by elements of $B_3$ described in Theorem \ref{one} is very similar to the picture obtained by physicists studying cascades of quiver gauge theories  on $X=\omega_{\PP^2}$ \cite{FHHI}. To obtain an exact match one should consider the subcategories $\A(g)\subset \D$ up to the action of the group $\Aut(\D)$. In the physicists' pictures the operation corresponding to tilting is called Seiberg duality \cite{BD,He} and the resulting webs are called duality trees. Physicists, particularly Hanany and collaborators, have computed many more examples (see for example \cite{BHK,FHHI,FHH}). 

\subsection{Resolutions of Kleinian singularities}
\label{dynkin}

There are a couple of interesting examples where the theory described above enable one to completely describe a connected component of
the space of stability conditions.

Let $G\subset\SL(2,\C)$ be a finite group and let $f\colon X\to Y$ be the minimal resolution of the corresponding
Kleinian singularity $Y=\C^2/G$. Define a full subcategory
\[\D=\{E\in\D^b\Coh(X):\mathbf{R} f_*(E)=0\}\subset\D^b\Coh(X).\]
The groups $G$ have an ADE classification so we may also consider the associated complex simple Lie algebra $\g=\g_{\C}$ with its Cartan subalgebra $\h\subset\g$ and root system $\Lambda\subset \h^*$.
The Grothendieck group $K(\D)$ with the Euler form can be identified with the root lattice $\Z\Lambda\subset \h^*$ equipped with the Killing form.

It was proved in \cite{Br4} (see also \cite{Th2}) that a connected component $\Stab^\dagger(\D)\subset \Stab(\D)$ is a covering space of
\[\hreg=\{\theta\in \h:\theta(\alpha)\neq 0\text{ for all }\alpha\in \Lambda\}.\]
The regions corresponding to stability conditions with a fixed heart are precisely the connected components of the inverse images of the
complexified Weyl chambers.
If we set $\Aut^\dagger(\D)$ to be the subgroup of $\Aut(\D)$ preserving this connected component one obtains
\[\frac{\Stab^\dagger(\D)}{\Aut^\dagger(\D)}\isom \frac{\hreg}{W^e},\]
where $W^e=W \rtimes\Aut(\Gamma)$ is the semi-direct product of the Weyl group of $\g$ with the finite group of automorphisms of the corresponding Dynkin graph.

In the same geometric situation  one can instead consider the full subcategory $\Dh\subset \D^b\Coh(X)$
consisting of objects supported on the exceptional locus of $f$.
A connected component of the space of stability conditions is then a covering space of the regular part of the affine
Cartan algebra $\hh$ and
\[\frac{\Stab^\dagger(\Dh)}{\Aut^\dagger(\Dh)}\isom \frac{\hregh}{\Wh^e},\]
where now $\Wh^e=\Wh \rtimes\Aut(\Gammah)$ is the semi-direct product of the affine Weyl group of $\g$ with the finite group of automorphisms of the corresponding affine Dynkin graph. For more details see \cite{Br4}.

In the $A_n$ case the spaces $\Stab(\D)$ and $\Stab(\Dh)$ are known to be connected and simply-connected \cite{Ish}.
A similar but more difficult example involving resolutions of three-dimensional singularities has been considered by Toda \cite{To} (see also \cite{Br6}). He has also considered
three-dimensional Calabi-Yau categories defined by considering a formal neighbourhood of a fibre of a K3 or elliptic fibration \cite{To2}.


\section{Geometric structures on spaces of stability conditions}

\label{font}

This section will be of a more speculative nature than the previous ones. I shall try to use ideas from mirror symmetry to make a few remarks about what geometric structures the space of stability conditions should carry.

\subsection{Stability conditions and the stringy K{\"a}hler moduli space}

Let $X$ be a simply-connected Calabi-Yau threefold and $\D=\D^b\Coh(X)$. In Section \ref{birdtwo} it was argued that one should expect an embedding of the stringy K{\"a}hler moduli space
$\M_K(X)$ in the double quotient \[\Aut(\D)\backslash\Stab(X)/\C.\] It is tempting to suggest that these two spaces should be identified. In fact it is easy to see using Theorem \ref{lasty} that this could never be the case. Put simply, the space $\Stab (X)$ is too big and too flat.

For concreteness let us take $X$ to be the quintic threefold. The stringy K{\"a}hler moduli space $\M_K(X)$ is, more or less by definition, the complex moduli space of the mirror threefold $Y$. As is well-known this is a twice-punctured two-sphere with a special point. The punctures are the large volume limit point and the conifold point, and the special point is called the Gepner point. The periods of the mirror $Y$ define holomorphic functions on $\M_{\C}(Y)$ which satisfy a third order Picard-Fuchs equation which has regular singular points at the special points.

Under mirror symmetry the periods of Lagrangian submanifolds of $Y$ correspond to central charges of objects of $\D$. Thus we see that the possible maps $Z\colon K(\D)\to\C$ occurring as central charges of stability conditions coming from points of $\M_K(X)$ satisfy the Picard-Fuchs equation for $Y$. Since these satisfy no linear relation, comparing with Theorem \ref{lasty} we see that the space $\Stab(X)$ must be four-dimensional 
 and the double quotient above is a three-dimensional space containing $\M_K(X)$ as a one-dimensional submanifold.
The embedding of this submanifold in $\Stab(X)$ is highly transcendental.

More generally, for a simply-connected Calabi-Yau threefold $X$ we would guess that the space $\Stab(X)$ is not the stringy K{\"a}hler moduli space, whose tangent space can be identified with $H^{1,1}(X)$, but rather some extended version of it, whose tangent space is 
\[\bigoplus_p H^{p,p}(X).\]
To pick out the K{\"a}hler moduli space as a submanifold of
\[\Aut(\D)\backslash\Stab(X)/\C.\]
we would need to define some extra structure on the space of stability conditions. Clues as to the nature of this extra structure can be obtained by studying other situations where extended moduli spaces occur in the mirror symmetry story.

\subsection{Extended moduli spaces}

There are (at least) three places where extended moduli spaces crop up in algebraic geometry: universal unfolding spaces, big quantum cohomology and the extended moduli spaces of Barannikov-Kontsevich.
All these spaces carry rich geometric structures closely related to Frobenius structures and all of them are closely related to moduli spaces of SCFTs. In each case one can make links with spaces of stability conditions, although none of these are close to being made precise. We content ourselves with giving the briefest outlines of the connections together with some references.

As explained by Takahashi \cite{Ta}, the unfolding space $T$ of an isolated hypersurface singularity $X_0$ of dimension $n$ should be related to the space of stability conditions on the Fukaya category of the Milnor fibre $X_t$ of the singularity. Note that
\[\mu=\dim_{\C} H_n(X_t,\C)=\dim_{\C} T.\]
Given a basis $L_1,\cdots, L_{\mu}$ of $H_n(X_t,\C)$,
K. Saito's theory of primitive forms shows that for a suitable family of holomorphic $n$-forms $\Omega_t$ on the fibres $X_t$ the periods
\[Z(L_i)=\int_{L_i} \Omega_t,\]
form a system of flat co-ordinates on the unfolding space $T$. Since these periods are the analogues of central charges this is exactly what one would expect from Theorem \ref{lasty}.

The big quantum cohomology of a Fano variety $Z$ seems to be related to the space of stability conditions on the derived categories of $Z$ and of the corresponding local Calabi-Yau variety $\omega_Z$. In the case when the quantum cohomology of $Z$ is generically semisimple Dubrovin showed how to analytically continue the prepotential from an open subset of $H^*(X,\C)$ to give a Frobenius structure on  a dense open subset of the configuration  space
\[M\subset \Con_n(\C)=\{(u_1,\cdots,u_n)\in\C^n:i\neq j\implies u_i\neq u_j\}.\]
According to a conjecture of Dubrovin \cite{Du2} the quantum cohomology of $Z$ is generically semisimple (so that one can define the above extended moduli space $M$) iff the derived category $\D^b\Coh(Z)$ has a full, strong exceptional
collection $(E_0,\cdots,E_{n-1})$ (so that one can understand the space of stability conditions by tilting as in Theorem \ref{one}). Moreover, in suitable co-ordinates, the Stokes matrix of the quantum cohomology $S_{ij}$ (which controls the analytic continuation of the Frobenius structure on $M$) is equal to the Gram matrix $\eu(E_i,E_j)$ (which controls the tilting or mutation process).  For more on this see \cite{Br5}.

Finally, Barannikov and Kontsevich \cite{BarKon} showed that if $X$ is a complex projective variety then the formal germ to deformations of $X$, whose tangent space has dimension $H^1(X,T_X)$, is contained in a larger formal germ whose tangent space has dimension \[\bigoplus_{p,q} H^p(X,\bigwedge^q T_X),\]
and which describes $A_\infty$ deformations of the category $\D^b\Coh(X)$.
Suppose $X_1$ and $X_2$ are a mirror pair of Calabi-Yau threefolds. Complex deformations of $X_1$ correspond to K{\"a}hler deformations of $X_2$. Passing to extended moduli spaces one might imagine that some global form of Barannikov and Kontsevich's space parameterising deformations of $\D^b\Coh(X_1)$ should be mirror to the space of stability conditions on $X_2$. Very schematically we might write
\[\Def(\D^b\Coh(X_1)) \isom \Stab(\D^b(\Coh(X_2)),\]
although to make the dimensions add up one should extend $\Stab(X_2)$ so that its tangent space is the whole cohomology of $X_2$ as in   Remark \ref{pc}.
Note that such an isomorphism would be a mirror symmetry statement staying entirely within the realm of algebraic geometry.

\subsection{An example}
Consider again the example of Section \ref{dynkin} relating to the resolution of the Kleinian singularity $Y=\C^2/G$. Consider the function on $\Stab(\D)$ defined by
\[\tag{$*$} F(Z)=\big.\sum_{\alpha\in\Lambda_+} Z(\alpha)^2 \log Z(\alpha)\]
Here the sum is over all positive roots and $Z$ denotes the central charge of a given point of $\Stab(\D)$. Of course, unless $\Stab(\D)$ is the universal cover of $\hreg$ the function $F$ may be many-valued. This function $F$ satisfies the WDVV equation (see for example \cite{Ve}). We can therefore use the derivatives of $F$ to define double and triple point functions
\[\langle\theta_1,\theta_2\rangle=\big.\sum_{\alpha\in\Lambda_+} \theta_1(\alpha)\theta_2(\alpha)\]
\[\langle\theta_1,\theta_2,\theta_3\rangle=\big.\sum_{\alpha\in\Lambda_+} \frac{\theta_1(\alpha)\theta_2(\alpha)\theta_3(\alpha)}{Z(\alpha)}\]
where $\theta_i\colon K(\D)\to\C$ are tangent vectors to $\Stab(\D)$ at a point $\sigma=(Z,\P)$.
Defining a multiplication on tangent vectors by
\[\langle\theta_1*\theta_2,\theta_3\rangle=\langle\theta_1,\theta_2,\theta_3\rangle=\langle\theta_1,\theta_2*\theta_3\rangle\]
defines an associative multiplication on the tangent bundle to $\Stab(\D)$ whose identity is the vector field $Z$. Since this identity is not flat this does not quite define a Frobenius manifold, but rather forms what Dubrovin calls an almost Frobenius manifold \cite{Du4}. In fact this structure is the almost-dual of the Frobenius manifold of Saito type on the unfolding space of the surface singularity $X$ (see \cite[Sections 5.1, 5.2]{Du4}).

It would be nice to generalise the function $F$ to some other examples.
In \cite{DGKV} there are formulae for prepotentials of gauge theories which look like $(*)$ (see for example Equation (3.1)) but with correction terms
involving sums over graphs. One of these was checked to satisfy the WDVV equation in \cite{CMMV}. This connection looks worthy of further investigation.

As a final remark, very recently Joyce \cite{Jo2} has constructed a flat connection on the space of stability conditions on an abelian category satisfying the Calabi-Yau condition. Extending this work to derived category seems to be problematic, but at present Joyce's approach seems to be our best hope for defining interesting structures on the space of stability conditions.

\subsection*{Acknowledgements}
The author is supported by a Royal Society University Research Fellowship.





\end{document}